\documentclass{daj}
\usepackage{enumitem,graphicx,xcolor,amssymb,amsmath,amsthm}

\newtheorem{thm}{Theorem}
\newtheorem{lem}[thm]{Lemma}
\newtheorem{cor}[thm]{Corollary}

\newcommand{\LARG}{\mathrm{LARG}}
\newcommand{\N}{\mathbb N}
\newcommand{\R}{\mathbb R}
\newcommand{\Z}{\mathbb Z}
\newcommand{\E}{\mathbb E}
\newcommand{\F}{\mathcal F}
\newcommand{\cross}{\mathsf{cross}}
\newcommand{\CD}{\mathsf{CD}}
\newcommand{\PP}{\mathbb P}
\newcommand{\Var}{\mathrm{Var}}

\newcommand{\cdinf}{\textsf{CD-inf}}

\dajAUTHORdetails{%
  title = {Geometric Random Graphs and Rado Sets of Continuous functions}, 
  author = {Anthony Bonato, Jeannette Janssen, and Anthony Quas},
  runningtitle = {Geometric Random Graphs in the Continuous Functions]},
  keywords = {geometric random graph, Rado graph},
}  

\dajEDITORdetails{%
   year={2021},
   number={3},
   received={2 March 2018},   
   revised={25 May 2020},    
   published={7 April 2021},  
   doi={10.19086/da.21473},       
}   

\begin{document}

\begin{frontmatter}[classification=text]

\title{Geometric Random Graphs and Rado Sets of Continuous functions} 

\author[ab]{Anthony Bonato\thanks{Supported by NSERC}}
\author[jj]{Jeannette Janssen\thanks{Supported by NSERC}}
\author[AQ]{Anthony Quas\thanks{Supported by NSERC}}

\begin{abstract}
We prove the existence of geometric Rado sets in the Banach space of continuous functions on $[0,1]$.
Starting from a countable dense set $S$, we form a random graph by independently joining pairs of elements of $S$
that are less than 1 apart by an edge with probability $p\in(0,1)$. If any two realizations of this
random graph are almost surely isomorphic, then the set $S$ is said to be \emph{geometrically Rado}.

We show that for suitable measures which we construct, almost all countable dense sets
of piecewise linear functions and of polynomials are geometrically Rado. The same holds for
countably many independent samples of integrated Brownian motion.
Moreover the almost sure random graph obtained in the three cases is the same.

We also study Brownian motion paths. Again, we show for a natural measure on Brownian paths,
almost all countable subsets are geometrically Rado
and the resulting graphs are almost surely of a unique isomorphism type.
However, the almost sure Rado graph in the Brownian motion case is not isomorphic
to the Rado graph in the previous cases.
\end{abstract}
\end{frontmatter}


\section{Overview}
We consider \emph{geometric graphs}, where the vertex set is a dense subset of a
separable metric space, and the edge set has the property that there exists $a>0$ such that
if the vertices $x$ and $y$ are joined by an edge, then $0<d(x,y)<a$.
In this paper the metric space will be a separable Banach space, $X$,
and $a$ will be 1.

If $V$ is a countable dense subset of $X$, and $p\in(0,1)$ is a parameter,
then we define a collection of random geometric graphs with vertex set $V$ by independently joining
each pair $x$ and $y$ in $V$ by an edge with probability $p$ if $0<\|x-y\|<1$.
We call this distribution the \emph{local area random graph} on $V$ with parameter $p$
and denote it by $\LARG(V,p)$.

Our main interest in this paper is in questions of isomorphism between random geometric graphs. First, we ask, for a fixed countable dense vertex set $V$, whether the graphs in $\LARG(V,p)$ are almost surely isomorphic. If they are, then we say that the set $V$ has the \emph{geometric Rado property}, by analogy with the discovery that random subgraphs on
a countable vertex set where all pairs of vertices are independently joined with probability $p$
are almost surely isomorphic to a universal graph.

Second, if a Banach space has countable dense subsets with the Rado property, we ask whether
there is a natural measure $\mu$ on $X$ such that for $\mu^\N$-almost every
countable dense $V$ (simply the set of elements appearing in the sequence in $X^\N$),
the elements of $\LARG(V,p)$ are almost all isomorphic. In this case, we also
ask whether for two independently chosen countable dense subsets $V$ and $V'$ and for all $p,p'\in(0,1)$, it is the case that almost every $G$ in $\LARG(V,p)$ and $G'\in\LARG(V',p')$ are isomorphic.
If it is, then we say that $\mu^\N$ is \emph{universally geometrically Rado}.
If $X$ supports natural measures generating universally geometrically Rado vertex sets, then it is
natural to ask if different measures lead to different universal graphs.

In \cite{BJ1} where these notions were introduced, Bonato and Janssen showed the surprising result that
if a finite-dimensional $\ell^\infty$ space is equipped with an absolutely continuous probability
measure $\mu$ with full support (that is, every open set has positive measure), then $\mu^\N$ is
universally geometrically Rado, and also the isomorphism class does not depend on the measure $\mu$.
The fact that $\mu$ is absolutely continuous guarantees that, almost surely, all of the points
sampled are \emph{integer distance free}, that is that no two distinct points have a coordinate
in which the points differ by
an integer. In fact, their proof shows that if $V$ and $W$ are dense countable integer distance free subsets
of $\ell^\infty$, then elements of $\LARG(V,p)$ and $\LARG(W,p')$ are almost surely isomorphic. 
They also show that in $\R^2$ equipped with the Euclidean norm, countable dense sets never
have the Rado property. 

In the elegant paper \cite{BBGLW}, Balister,
Bollob\'as, Gunderson, Leader, and Walters
give a necessary and sufficient criterion for the existence of Rado sets in finite-dimensional Banach spaces.
They uniquely express $X$ as an ``$\infty$-decomposition" $(U\oplus\ell^d_\infty)_\infty$ into a
maximal $\ell^\infty$ part and the rest. They then show that the space
admits geometric Rado sets if and only if $d>0$, and it has the universal geometric Rado property if and only if $U=\{0\}$
(that is, if and only if $X$ is isometrically
isomorphic to $\ell^d_\infty$ for some $d$).
They also show that there is a countable dense subset $V$ of $\ell^\infty$ that is not integer distance
free which does not have the geometric Rado property, thereby emphasizing the role of
the integer distance free property. 

In \cite{BJQ1}, we studied the existence of Rado sets in the separable infinite-dimensional
subspaces of $\ell^\infty$: $c_0$ (of sequences converging to 0) and $c$ (of convergent sequences).
We exhibited fully supported measures on $c_0$ and $c$  with the universal geometric Rado property
and showed that in both cases, the universal graph obtained was the same for all measures in a
reasonably general class (but the universal graphs for $c_0$ and $c$ are not isomorphic; see
Corollary~\ref{cor:noniso}). In \cite{AS}, the authors considered Rado sets in metric spaces on circles of a fixed
circumference.

In this paper, we study the Banach space $C[0,1]$ of continuous real-valued functions on $[0,1]$,
a separable subspace of $L^\infty[0,1]$.
We exhibit probability measures $\mu$ on $C[0,1]$ supported on polynomials,
on piecewise linear functions, and on integrated Brownian motions,
such that $\mu^\N$ is universally geometrically Rado. Further, the universal random graph is the same for
these classes. Additionally, we exhibit a probability measure $\nu$ on $C[0,1]$ supported on Brownian motion
functions such that $\nu^\N$ is universally geometrically Rado.
However, it turns out that the universal random
graph for the Brownian motion functions is not isomorphic to the
universal random graph for the polynomials and piecewise
linear functions.
Informally, the universal random graphs ``see the intersections between pairs of functions". The fact that
pairs of Brownian motions cross infinitely often if they cross at all leads to non-isomorphism of the graphs.

The current paper revolves around two themes:
\begin{enumerate}
\item
showing that the $\ell^\infty$-like space $C[0,1]$ supports natural measures on countable dense subsets
that exhibit the universal geometric Rado property;
\item
understanding the geometric information about countable dense subsets that may be obtained
from the isomorphism type of the local area random graphs that they support (in particular
exhibiting two natural classes of distributions on $C[0,1]$ with the universal geometric
Rado property, and showing that their universal Rado graphs fail to be isomorphic).
\end{enumerate}

The \emph{crossings} between elements functions $f$ and $g$ (that is, points
$x\in[0,1]$ such that $f(x)-g(x)\in\Z$) play a key role in all of this that is similar to the
role played by the \emph{integer distance free} property in \cite{BJ1}.

\section{The Banach space and random graph geometries}

A geometric graph $G$ is said to be
\emph{well-connected} if it has the property that for every finite subset $S$ of $V(G)$ and
for every non-empty open set $U\subseteq X$ such that $\|u-x\|<1$ for all $u\in U$ and $x\in S$,
there exists $y\in U\cap V(G)$ such that $y$ is adjacent to each $x\in S$.
By a Borel-Cantelli argument,
it is straightforward to show that if $V$ is a countable  dense subset of $X$,
then almost every $G$ sampled from $\LARG(V,p)$
is well-connected.

\begin{lem}\label{lem:wellcon}
If $G$ and $G'$ are well-connected geometric
graphs in $X$ and $X'$ and $\theta\colon G\to G'$ is a graph isomorphism, then we have the following:
\begin{enumerate}
\item $\lfloor \|x-y\|\rfloor=\lfloor \|\theta(x)-\theta(y)\|\rfloor$ for all $x,y\in V(G)$
satisfying $\|x-y\|\ge 2$; and
\item $\lceil \|x-y\|\rceil=\lceil \|\theta(x)-\theta(y)\|\rceil$ for all $x,y\in V(G)$.
\end{enumerate}
In particular, if $\|x-y\|$ is an integer and is at least 2, then $\|\theta(x)-\theta(y)\|=\|x-y\|$.
\end{lem}

For the time being, in the infinite-dimensional case, we are unable to rule out
the possibility that a pair of vertices
separated by a distance of 1 are mapped to vertices separated by a distance less than 1
that are not connected by an edge. However, this can be ruled out in the finite-dimensional case; see
Corollary~\ref{cor:fdim} below.

The following corollary of Lemma~\ref{lem:wellcon} was established in \cite{BJQ1}. The corollary 
was established using a theorem of Dilworth \cite{Dilworth} showing that approximately surjective approximate isometries
may be uniformly approximated by surjective isometries, and the theorem of Mazur and Ulam stating that
surjective isometries of real Banach spaces may be expressed as the composition of a translation
and a linear isometry.

\begin{cor}\label{cor:noniso}
If $G$ and $G'$ are isomorphic well-connected geometric graphs in $X$ and $X'$, then the Banach spaces $X$
and $X'$ are isometrically isomorphic.
\end{cor}

The corollary implies that a well-connected geometric graph with vertices in a separable Banach space
determines the Banach space structure.

A map between two sets $V$ and $V'$ is called a \emph{step-isometry} if
$\lfloor \|\theta(x)-\theta(y)\|\rfloor=
\lfloor\|x-y\|\rfloor$ for all $x,y\in V$. As shown in Lemma~\ref{lem:wellcon}, a graph isomorphism
between two well-connected graphs is almost a step-isometry (modulo the possibility mentioned above).
We use this in the context of $C[0,1]$
to build graph isomorphisms: given graphs $G$ and $G'$ sampled from $\LARG(V,p)$ and
$\LARG(V',p')$, we build a graph isomorphism using a back-and-forth technique, ensuring that at each
stage of the construction (at which time only finitely many vertices are matched), the part of the
isomorphism so far constructed is a step-isometry. As long as all of the vertices in both graphs are exhausted
by this process, this yields an isomorphism between $G$ and $G'$.

The above shows (as first observed in [BJ1]),
there is an intimate connection between
step-isometries and geometric graph isomorphisms. This is elaborated upon in \cite{BBGLW}; in particular, 
it is shown that if a finite-dimensional space $X$ has $\infty$-decomposition
$(U\oplus \ell_\infty^d)_\infty$, then a step-isometry may be factorized over the
decomposition as a surjective linear isometry (on $U$) and a step-isometry on the $\ell^\infty$
factor. Establishing the Rado property depends on the availability of a plentiful supply of
step-isometries and this gives rise to the failure of the Rado property on
spaces that are not isomorphic to $\ell_\infty^d$.

\begin{proof}[Proof of Lemma~\ref{lem:wellcon}]
It was shown in \cite{BJ1} that if $x,y$ are elements of $V(G)$ for any integer $k\ge 2$, $\|x-y\|<k$
implies $d_G(x,y)\le k$. To see this, consider the $k-1$ equally spaced points
$x+\frac{i}{k}(y-x)$ for $i=1,\ldots, k-1$ along the chord $[x,y]$.
Since $\frac{\|y-x\|}{k}<1$, there are open balls around these points
with the property that each point in the $i$th ball is less than one unit distant from each point in the next ball
and the first and last balls are within 1 unit of $x$ and $y,$ respectively.
For $i=1$, by well-connectedness there exists a vertex $x_1$ in the first ball adjacent to $x$.
Continuing inductively one may select $x_1,\ldots,x_{k-2}$ such that each is adjacent to the next.
Finally, $x_{k-1}$ is chosen to be adjacent to both $x_{k-2}$ and $y$.
On the other hand, it is straightforward that for any $k\in\N$,
$d_G(x,y)\le k$ implies $\|x-y\|<k$.
It follows that for any $k\ge 3$, $d_G(x,y)=k$ if and only if $\|x-y\|\in [k-1,k)$.
In particular, if $\|x-y\|\ge 2$, then $\lfloor \|x-y\|\rfloor=d_G(x,y)-1=d_{G'}(\theta(x),\theta(y))-
1=\lfloor\|\theta(x)-\theta(y)\|\rfloor$.

For the second part, suppose $\|x-y\|>k\in\N_0$. By well-connectedness,
there exist $w$ adjacent to $x$
and $z$ adjacent to $y$ such that $\|w-z\|>k+2$
(for example, let $u=(y-x)/\|y-x\|$ be the direction joining $x$ to $y$
and choose $w$ in a small ball around $x-u$ and $z$ in a small ball around $y+u$).
Hence, we have that $d_G(w,z)\ge k+3$, so that
$d_{G'}(\theta(w),\theta(z))\ge k+3$ and
by the previous part, $\|\theta(w)-\theta(z)\|\ge k+2$.
Now, it follows that $\|\theta(x)-\theta(y)\|\ge \|\theta(w)-\theta(z)\|-\|\theta(w)-\theta(x)\|
-\|\theta(z)-\theta(y)\|>k$.
\end{proof}

Recall that a sequence $(y_n)$ in a Banach space $X$ is said to \emph{converge weakly} to $y$
if $\phi(y_n)\to \phi(y)$ for all $\phi$ in $X^*$, the dual space of $X$.

\begin{lem}\label{lem:contin}
Let $G$ and $G'$ be well-connected graphs on the Banach space $X$ and suppose that
they are isomorphic by the map $\theta$.
If $(x_n)$ is a sequence of elements of $V(G)$ such that $x_n\to x\in V(G)$,
then $\theta(x_n)$ converges weakly to $\theta(x)$ in $V(G')$.

In particular, if $X$ is finite-dimensional, $\theta$ is continuous.
\end{lem}

\begin{proof}
We prove this by contradiction. Let $y=\theta(x)$ and $y_n=\theta(x_n)$.
Suppose that $\phi$ is a linear functional of norm 1 so that
$\phi(y_n)\not\to\phi(y)$. Since $x_n\to x$, all but finitely many $x_n$'s lie in the closed ball of
radius 1 about $x$, so that by the second statement of Lemma~\ref{lem:wellcon},
for all except these finitely many $n$, $y_n$ lies in the closed ball of radius
1 about $y$. Hence, we may extract
a subsequence such that $\phi(y_n)\to \phi(y)-a$, where we may assume,
by negating $\phi$ if necessary, that
$a\in(0,1)$. Since $\|\phi\|=1$ and $V(G')$ is dense, there exists
$z=\theta(w)\in V(G')$ such that $\|z-y\|<2$; $\phi(z)-\phi(y)>2-\frac a2$. For
all sufficiently large $n$, we have that $\phi(z)-\phi(y_n)>2$, so that $\|z-y_n\|> 2$. By the
second statement of Lemma~\ref{lem:wellcon}, $\|w-x_n\|>2$,
so that $\|w-x\|\ge 2$. We then have that $d_G(w,x)\ge 3$. However, $\|z-y\|<2$ implies
that $d_{G'}(z,y)=2$ by the first paragraph of the proof of Lemma~\ref{lem:wellcon},
which contradicts the assumption that $G$ and $G'$ are isomorphic.

If $X$ is finite-dimensional, then weak convergence is equivalent to convergence in norm.
\end{proof}

\begin{cor}\label{cor:fdim}
If $X$ is a finite-dimensional Banach space and $G$ and $G'$ are two well-connected
geometric graphs on $X$, that are isomorphic by a map $\theta$, then $\|\theta(x)-\theta(y)\|=1$
if and only if $\|x-y\|=1$.

Hence, in the finite-dimensional case,
for any $k\in\N_0$, $\|x-y\|=k$ if and only if $\|\theta(x)-\theta(y)\|=k$;
$k<\|x-y\|<k+1$ if and only $k<\|\theta(x)-\theta(y)\|<k+1$.
\end{cor}

\begin{proof}
Let $x,y\in V(G)$ satisfy $\|x-y\|=1$. By Lemma~\ref{lem:wellcon}, we know that 
$\|\theta(x)-\theta(y)\|\in(0,1]$. Let $x_n\to x$ and $y_n\to y$ with $x_n,y_n\in V(G)$
and $\|y_n-x_n\|>1$ (the existence of such sequences is guaranteed by density of $V(G)$:
let $u=y-x$ be the unit vector joining $x$ and $y$;
let $x_n$ be a vertex in a ball of radius $\frac 1{2n}$ around
$x-\frac un$ and let $y_n$ be a vertex in a ball of radius $\frac 1{2n}$ around
$y+\frac un$). We then have that $\|\theta(y_n)-\theta(x_n)\|>1$
by the second statement of Lemma~\ref{lem:wellcon}
and $\theta(x_n)\to\theta(x)$ and $\theta(y_n)\to\theta(y)$ by
Lemma~\ref{lem:contin}, so that $\|\theta(y)-\theta(x)\|\ge1$,
completing the proof of the first part. The second part follows by combining this with Lemma~\ref{lem:wellcon}.
\end{proof}

We know from Lemma~\ref{lem:wellcon} that the isometries of well-connected
geometric graphs give rise to step-isometries of the vertex sets. It is therefore
useful to have a large collection of step-isometries of $C[0,1]$. The following lemma gives a class of step-isometries that underlie the inductive
proof that we give of the Rado properties.

Informally, we can think of the step-isometries that we build as acting on graphs of functions on
$[0,1]$ as follows. We first apply a homeomorphism to the $x$ coordinates.
We then apply a step-isometry (depending continuously on $x$) to the graphs
of the functions over the $x$ coordinate.
This should be compared with the Banach-Stone theorem; see Theorem~\ref{thm:BS}.

\begin{lem}\label{lem:stepisos}
Let $\Phi\colon [0,1]\times\R\to \R$ be a continuous function such that
$\Phi(y,t)$ is a strictly increasing function of $t$ satisfying $\Phi(y,t+1)=\Phi(y,t)+1$ for all $y\in [0,1]$,
$t\in\R$.
If $\psi\colon[0,1]\to[0,1]$ is an increasing continuous bijection,
then the map $\Theta\colon C[0,1]\to C[0,1]$ defined by
$\Theta[f](y)=\Phi(y,f(\psi(y)))$ is a bijective continuous step-isometry of $C[0,1]$.
\end{lem}

\begin{proof}
We first establish that $\Theta$ is a homeomorphism of $C[0,1]$, and that its inverse is of the
same form.
Note that the map $f\mapsto f\circ\psi$ is a bijective linear isometry
of $C[0,1]$ and that $\Theta$ may be expressed as a composition of $g\mapsto \Phi(y,g(y))$
with this bijection. Hence, it suffices to consider the case where $\psi(y)=y$ and $\Theta[f](y)=
\Phi(y,f(y))$.
Let $\phi_y$ be defined by $\phi_y(t)=\Phi(y,t)$. By hypothesis, $\phi_y$ is a
bijection of $\R$.

The inverse of $\Theta$ is given by $\Theta^{-1}(g)(y)=
\phi_y^{-1}(g(y))$. We claim that $\Theta^{-1}$ is continuous.
To see this, let $\epsilon>0$ and let $\delta=\min_{y\in[0,1]}\min_{t\in\R}[\phi_y(t+\epsilon)-\phi_y(t)]$
(this quantity is positive by the fact that $\phi_y$ is periodic-plus-identity and by compactness of [0,1]).
Now for any $y\in [0,1]$ and $a\in\R$,
$\phi_y(\phi_y^{-1}a-\epsilon)\le \phi_y(\phi_y^{-1}a)-\delta=a-\delta$, so that
$\phi_y^{-1}(a-\delta)\ge\phi_y^{-1}a-\epsilon$ and similarly
$\phi_y^{-1}(a+\delta)\le\phi_y^{-1}a+\epsilon$.
In particular, if $|a-b|<\delta$, then $|\phi_y^{-1}b-\phi_y^{-1}a|<\epsilon$.
Now if $\|f-g\|<\delta$, we see that $\|\Theta^{-1}f-\Theta^{-1}g\|<\epsilon$. (We have shown
in this paragraph that the inversion map from the set of bijective
homeomorphisms of the line such that $\phi(t+1)=\phi(t)+1$ to itself is continuous.)

If $\lfloor \|f-g\|\rfloor=k$, then we have that $|f(x)-g(x)|<k+1$ for each $x$.
We note that $\Theta[f](y)-\Theta(g)(y)=\Phi(y,f(y))-
\Phi(y,g(y))$. Since $-(k+1)<f(y)-g(y)<k+1$,
it follows that $-(k+1)<\Phi(y,f(y))-\Phi(y,g(y))<k+1$ for each $x$,
so that $\lfloor \|\Theta(f)-\Theta(g)\|\rfloor\le k=\lfloor \|f-g\|\rfloor$. Applying this argument to
$\Theta^{-1}$ (which we have already shown to be of the same form),
we deduce that $\lfloor\|\Theta(f)-\Theta(g)\|\rfloor=\lfloor\|f-g\|\rfloor,$
as required.
\end{proof}

\section{Transversely dense functions}
Lemma~\ref{lem:wellcon} suggests that pairs of functions separated by an
integer (norm) distance are particularly problematic.
In fact more is true: not only the maximum distance is important,
but also places where pairs of functions differ by an integer (even if that integer is
not the largest distance between the functions).
A key issue underlying most of this paper is that any graph isomorphism
has to preserve ``crossing patterns". Two continuous functions $f$ and $g$
\emph{cross} at $x$ if $f(x)-g(x)\in\Z$. We define $\cross(f,g)$ to be
$\{x\colon f(x)-g(x)\in\Z\}$. If $\mathcal F$ is a set of functions, then $\cross(\mathcal F)$ denotes $\bigcup_{f,g\in\mathcal F; f\ne g}
\cross(f,g)$. If $\mathcal F$ is a collection of continuous functions, then
notice that for all $f,g\in\mathcal F$, $\lfloor f-g\rfloor$ is constant on each
component of $[0,1]\setminus\cross(\mathcal F)$.

For the first class of distributions on $X^\N$ that we show to be geometrically Rado,
we impose the condition that crossing sets are discrete and disjoint.
A subset $V$ of $C[0,1]$ is said to be \emph{transverse} if the following two properties hold:
\begin{enumerate}[label=(Tr\arabic*)]
\item For each pair of distinct functions, $f,g\in V$, $\cross(f,g)$ is finite;
does not contain 0 or 1; and does not
contain any local maxima or minima of $f-g$.
\label{it:Tr1}
\item The sets $\cross(f,g)$ are disjoint as $(f,g)$ runs over the pairs of distinct functions in $V$.
\label{it:Tr2}
\end{enumerate}

Given a finite transverse collection of continuous functions, $f_1,\ldots,f_m$,
its \emph{crossing data}, $\CD(f_1,\ldots,f_m)$ consists of
a \emph{crossing tuple} $(x_1,\ldots,x_{k-1})$ and an array
$(d_{ab}^{(i)})_{1\le a<b\le m;\ 1\le i\le k}$ of integers
so that $$\cross(\{f_1,\ldots,f_m\})=\{x_1,\ldots,x_{k-1}\};$$ $x_i<x_{i+1}$ for each $i$; and
$(f_b-f_a)$ takes values in the range $(d_{ab}^{(i)},d_{ab}^{(i)}+1)$ for
$x\in (x_{i-1},x_i)$ where $x_0$ is taken to be
0 and $x_{k}$ is taken to be 1.

If $\CD(f_1,\ldots,f_m)=\Big((x_1,\ldots,x_{k-1}),(d_{ab}^{(i)})\Big)$ and $\CD(g_1,\ldots,g_m)=
\Big((x_1',\ldots,x'_{k'-1}),({d'}_{ab}^{(i)})\Big)$, then we write
$\CD(f_1,\ldots,f_m)\sim_\epsilon \CD(g_1,\ldots,g_m)$ if $k=k'$;
$|x_i-x'_i|<\epsilon$ for each $i$; and ${d'}_{ab}^{(i)}=d_{ab}^{(i)}$ for all $1\le a<b\le m$,
$1\le i\le k$. If $k=k'$ and ${d'}_{ab}^{(i)}=d_{ab}^{(i)}$ for all $i,a,b$ with no restriction on $|x_i-x'_i|$, then we
write $\CD(f_1,\ldots,f_m)\sim_1 \CD(g_1,\ldots,g_m)$.

Finally, a countable subset, $V$, of $C[0,1]$ is \emph{transversely dense} if the following properties hold:
\begin{enumerate}[label=(TD\arabic*)]
\item $V$ is transverse.
\item \label{it:TD2}
For each finite collection $\mathcal F=\{f_1,\ldots,f_m\}$ in $V$, for
each $f\in C[0,1]\setminus \mathcal F$ such that
$\mathcal F\cup\{f\}$ is transverse, for each $\epsilon>0$, there exist
infinitely many $g\in V\setminus\mathcal F$
such that $\|f-g\|<\epsilon$ and $\CD(f_1,\ldots,f_m,g)\sim_\epsilon \CD(f_1,\ldots,f_m,f)$.
\end{enumerate}
In particular, a transversely dense subset of $C[0,1]$ is dense.
We say a sequence $(f_1,f_2,\ldots)$ of distinct functions is transversely dense if
$\{f_1,f_2,\ldots\}$ is transversely dense.
A function $g$ satisfying \ref{it:TD2} is called an $\epsilon$-\emph{approximation} of $f$ relative to $\mathcal F$.
These $\epsilon$-approximations are the key to the back-and-forth technique that we use.

\begin{thm}\label{thm:TD}
Let $V$ and $W$ be countable transversely dense subsets of $C[0,1]$ and let $0<p,p'<1$.
If $G$ and $G'$ are sampled from $\LARG(V,p)$ and $\LARG(W,p'),$ respectively,
then $G$ and $G'$ are almost surely isomorphic.

In particular, there exists a graph $\mathsf G_{TD}$ such that
if $\mu$ is a measure on $X$ where $\mu^\N$-a.e.\ sequence is
transversely dense, then $\mu$ is universally geometrically Rado
and for $\mu^\N$-a.e.\ vertex set $V$, almost every element of
$\LARG(V,p)$ is isomorphic to $\mathsf G_{TD}$ for any $0<p<1$.
\end{thm}

\begin{proof}
Let $V=\{f_1,f_2,\ldots\}$ and $W=\{g_1,g_2,\ldots\}$ and let $G$ and $G'$ be sampled from
$\LARG(V,p)$ and $\LARG(W,p'),$ respectively. We describe a
countable-step algorithm to build an isomorphism
from $G$ to $G'$ that fails at each time-step with probability 0; and if it does not fail at any time-step,
yields an isomorphism of $G$ and $G'$.

We use a back-and-forth process to (almost surely) build two sequences
$(i_n)$ and $(j_n)$ with the properties:
\begin{enumerate}[label=(T\ref{thm:TD}-\arabic*)]
\item\label{it:Ind1}
For each $n\in\N$, the graphs induced by $G$ on $\{f_{i_1},\ldots,f_{i_n}\}$ and
by $G'$ on $\{g_{j_1},\ldots,g_{j_n}\}$ are isomorphic by the map sending $f_{i_l}$ to $g_{j_l}$;
\item \label{it:Ind2}For each $n\in \N$, $\{i_1,\ldots,i_n\}$ and
$\{j_1,\ldots,j_n\}$ contain $\{m\in\N\colon m\le \frac n2\}$;
\item \label{it:Ind3}$\CD(f_{i_1},\ldots,f_{i_n})\sim_1 \CD(g_{j_1},\ldots,g_{j_n})$.
\end{enumerate}
It is evident that if properties \ref{it:Ind1} and \ref{it:Ind2} are satisfied, then $G$ and $G'$ are isomorphic.
Property \ref{it:Ind3} ensures that the partial isomorphism can always be extended.

To describe the algorithm, we initially set up two queues: $Q$ and $Q'$,
initially consisting of $(f_1,f_2,\ldots)$
and $(g_1,g_2,\ldots)$. As usual in the back-and-forth process, we take an unmatched
vertex in one of the graphs and match it to a suitable unmatched vertex in the other graph. To find
a suitable vertex, transverse density is shown to yield infinitely many suitable vertices to match to (based
only on the crossing pattern). We then apply
the Borel-Cantelli lemma to ensure that we pick a vertex with the correct adjacencies to the
previously matched vertices to ensure that the enlarged matching remains a graph isomorphism.
A critical hypothesis here is that all of the
events being tested in a single application of the lemma are independent, and it is helpful if they
all have equal probabilities. For this reason, when searching for a match for a particular vertex
in one of the graphs, we \emph{only test vertices in the other graph that have not previously been tested}.
That is, we take elements from the appropriate queue and stop when we find a match.
Elements that are tested are removed from the queue.
We make a further induction hypothesis.
\begin{enumerate}[resume*]
\item The queues $Q$ and $Q'$ consist of all but finitely many vertices of $G$ and $G'$.
\end{enumerate}

For the base case, we set $i_1=j_1=1$ and remove $f_1$ and $g_1$ from $Q$ and $Q'$.
We then inductively assume that
$i_1,\ldots, i_{n-1}$ and $j_1,\ldots,j_{n-1}$ have been chosen so that
the induction hypotheses are satisfied.
If $n$ is even, then let $i_n$ be the first unmatched vertex in $G$,
while if $n$ is odd, let $j_n$ be the first unmatched vertex
in $G'$. We then seek suitable $j_n$ or $i_n,$ respectively. Since the processes are essentially
identical, we deal only with the case where
$n$ is even, so we are looking for a suitable $j_n$.

We first identify a target function $g$ in $C[0,1]$, and then find a suitable approximation in $W$.
Let $\cross(f_{i_1},\ldots,f_{i_{n-1}})=\{x_1,\ldots,x_{k-1}\}$.
By the induction hypothesis \ref{it:Ind3}, we may
write  $$\cross(g_{j_1},\ldots,g_{j_{n-1}})=\{y_1,\ldots,y_{k-1}\},$$
and define $x_0=y_0=0$ and $x_k=y_k=1$.
We build the target function $g$ as follows. Let $\psi\colon[0,1]\to[0,1]$
be the continuous increasing piecewise linear map
sending each $[y_{l-1},y_l]$ to $[x_{l-1},x_l]$.
By \ref{it:Ind3}, we have $g_{j_a}(y)+p<g_{j_b}(y)+q$ if and only if $f_{i_a}(\psi(y))+p<f_{i_b}(\psi(y))+q$
for $a,b\in \{1,\ldots,n-1\}$ and $p,q\in\Z$.
Let $\phi_y$ be the piecewise linear map from $\R$ to $\R$ satisfying
\begin{equation}\label{eq:phidef}
\phi_y(f_{i_m}(\psi(y))+p)=g_{j_m}(y)+p
\end{equation}
for each $1\le m< n$ and $p\in\Z$ that is linear between consecutive
values of $f_{i_m}(y)+p$'s. This satisfies
$\phi_y(t+1)=\phi_y(t)+1$ for each $t\in\R$ and $y\in[0,1]$.
We write $f_{i,p}(x)$ for $f_i(x)+p$ and similarly $g_{j,p}(y)=g_j(y)+p$
for any $x,y\in[0,1]$ and $p\in\Z$.

The map obtained sending $f_{i_m,p}$ to $g_{j_m,p}$
is of the type described in Lemma~\ref{lem:stepisos} (also
described informally prior to the statement of that lemma).
The map $\psi$ is the piecewise linear homeomorphism sending the original
crossing points of the (finite collection of) $f$'s to the crossing points
of the corresponding $g$'s. The step-isometries over each $x$ coordinate
are piecewise linear, and send integer translates of each $f$ to the corresponding
integer translate of a $g$.

The fact that $\{f_{i_1},\ldots,f_{i_{n-1}}\}$ and
$\{g_{j_1},\ldots,g_{j_{n-1}}\}$ are transverse (specifically the fact
that $f_{i_a}$ and $f_{i_b}$ cross at $\psi(y)$ if and only if
$g_{j_a}$ and $g_{j_b}$ cross at $y$) ensures that the map $(y,t)\mapsto \phi_y(t)$
is continuous. To see this, notice that by \eqref{eq:phidef},
for $t\not\in S(y):=\{f_{i_a}(\psi(y))+p\colon 1\le a<n;\ p\in\Z\}$, the piecewise
linear function $\phi_y$ satisfies
$$
\phi_y(t)=\frac{[f_{i_M,P}(\psi(y))-t]g_{j_m,p}(y)}{f_{i_M,P}(\psi(y))-f_{i_m,p}(\psi(y))}
+
\frac{[t-f_{i_m,p}(\psi(y))]g_{j_M,P}(y)}{f_{i_M,P}(\psi(y))-f_{i_m,p}(y)},
$$
where $m$ and $p$ are such that $f_{i_m,p}(\psi(y))=\max(S(y)\cap (-\infty,t))$;
and $M$ and $P$ are such that
$f_{i_M,P}(\psi(y))=\min(S(y)\cap (t,\infty))$. If $t$ is close to $f_{i_a,b}(\psi(y_0))$
and $y$ is close to $y_0$,
then $\phi_y(t)$ is close to $g_{j_a,b}(y_0)$.

The target function is defined by
$$
g(y)=\phi_y(f_{i_n}(\psi(y))).
$$

By construction, the places where $g$ crosses $g_{j_1},\ldots,g_{j_{n-1}}$ are the pre-images under $\psi$
of the places where $f_{i_n}$ crosses $f_{i_1},\ldots,f_{i_{n-1}}$.
We see that $g(y)<g_{j_m,p}(y)$ if and only if $f_{i_n}(\psi(y))<f_{i_m,p}(\psi(y))$.
It follows that $\CD(f_{i_1},\ldots,f_{i_{n}})\sim_1\CD(g_{j_1},\ldots,g_{j_{n-1}},g)$.
Hence, by the transverse density hypothesis, there exist infinitely many $\ell$'s such that
$\CD(f_{i_1},\ldots,f_{i_{n}})\sim_1\CD(g_{j_1},\ldots,g_{j_{n-1}},g_\ell)$.
Notice that $\|f_{i_a}-f_{i_n}\|<1$ if and only if $\|g_{j_a}-g_\ell\|<1$ for each of these $\ell$'s and
each $1\le a\le n-1$.
Sequentially testing only those $\ell$'s (all but finitely many of these) that remain in $Q'$,
each has probability at least $\min(p',1-p')^{n-1}$ (where $p'$ is the parameter in $\LARG(W,p')$)
to have the same pattern of
adjacencies to $g_{j_1},\ldots,g_{j_{n-1}}$
as $f_{i_n}$ has to $f_{i_1},\ldots,f_{i_{n-1}}$, establishing the induction hypothesis
for $n$. The (almost surely finitely many)
tested $g$'s should then be removed from $Q'$ and $f_{i_n}$ should be removed from $Q$.
\end{proof}

We now exhibit two examples of measures on $C[0,1]$ satisfying the
hypotheses of the theorem, supported on polynomials
and piecewise linear functions, respectively. Let $P$ be the map from $\N\times \R^{\N_0}$ to $C[0,1]$ defined by $P(k,a_0,a_1,\ldots)=f$,
where $f(t)=\sum_{i=0}^k a_it^i$. Let $\N\times\R^{\N_0}$ be equipped with the product measure,
where the first coordinate is
geometrically distributed with parameter $\frac 12$ (that is, $\PP(N=k)=1/2^k$ for $k=1,2,\ldots$). The remaining coordinates are standard normal random variables.
Let $\mu_\text{poly}$ be the measure on
$C[0,1]$ obtained as the push-forward of $\PP$ under $P$ (that is, a randomly sampled element
is a polynomial of degree $k$ where $k$ is a geometric random variable whose coefficients are
sampled from independent standard normal random variables).

\begin{thm}\label{thm:polytd}
Let $\mu_\text{poly}$ be as above. For $\mu_\text{poly}^\N$-a.e.\ subset $V$ of $C[0,1]$,
and each $0<p<1$, almost every realization of $\LARG(V,p)$ is isomorphic to $\mathsf{G}_{TD}$.
\end{thm}

\begin{proof}
Let a randomly sampled $V$ be $\{f_1,f_2,\ldots\}$. We show that $V$ is
transversely dense, so that the desired
result follows from Theorem~\ref{thm:TD}

We first establish that $V$ is almost surely transverse.
Since $\mu_\text{poly}$ is a continuous measure, the $(f_i)$ are almost surely pairwise distinct, so that
$\cross(f_i,f_j)$ is finite since the functions are polynomials. If $f$ is sampled from $\mu_\text{poly}$,
we can consider $f(0)$ and $f(1)$ as random variables.
These also have continuous distributions, so the probability
that $\cross(f_i,f_j)$ intersects $\{0,1\}$ is 0 for each of the countably many $(i,j)$ pairs.

For the local maxima and minima of $f_i-f_j$, the (finitely many)
locations are determined by the non-constant terms of the
pair of $f_i(t)-f_j(t)$. Conditioned on the non-constant terms, the value of the difference at the locations of
the critical points is again continuously distributed, so that almost surely, $f_i-f_j$ does not
differ by an integer at any critical point for each of the countably many $(i,j)$ pairs.
We have established \ref{it:Tr1}.

Given a pair of functions $f_i$ and $f_j$, there are finitely many crossings. Now given the coefficients of the
polynomial $f_k$ except for the constant term, we see that conditioned
on all of this information, the probability
that $f_i$, $f_j$ and $f_k$ have a common crossing is 0.
A similar argument shows that the probability that two
other functions $f_k$ and $f_l$
cross at one of the crossings of $f_i$ and $f_j$ is 0, establishing \ref{it:Tr2}.

Finally, we establish \ref{it:TD2}.
Let $\mathcal F=\{f_1,\ldots,f_m\}$ (this may be done without loss of generality by expanding the finite
collection if necessary). Let $f\in C[0,1]$ be such that $\mathcal F\cup\{f\}$
 is transverse and let $\epsilon>0$.
We may assume without loss of generality that $\epsilon$ is less than $\frac 12\min\{|x-y|\colon
x,y\in \cross(f_1,\ldots,f_m,f),\ x\ne y\}$.
First, by the transverse assumption, we may replace $f$ by a
continuously differentiable function $g$ such that
$\CD(f_1,\ldots,f_m,g)=\CD(f_1,\ldots,f_m,f)$ and $(g-f_i)'$ is non-zero at each crossing of $g$ and $f_i$.
Applying the Weierstrass approximation theorem to $g'$,
there is a polynomial $q$, of degree $d$ say, and a $\delta>0$
such that if $\|Dh-q\|<\delta$ and $|h(0)-g(0)|<\delta$, then
$\|h-g\|<\epsilon$ and $\CD(f_1,\ldots,f_m,h)\sim_\epsilon
\CD(f_1,\ldots,f_m,g)$. Now there is an open collection of polynomials
of degree $d+1$ satisfying these constraints,
so that there are, as required, infinitely many $f_j$ such that $\CD(f_1,\ldots,f_m,f_j)\sim_\epsilon
\CD(f_1,\ldots,f_m,f)$.
\end{proof}

We now prove a similar result for piecewise linear functions.
Given an $n\in\N$, a sequence of distinct $(x_i)_{i\ge 1}$ taking
values in $(0,1)$  and a sequence of $(y_i)_{i\ge 0}$ taking
real values, let $\Phi(n,(x_i),(y_i))$ be the piecewise linear function with $n$ pieces joining
$(\hat x_{i-1},y_{i-1})$ to $(\hat x_i,y_i)$ for $i=1,\ldots,n$ by a straight line
where $\hat x_0=0$, $\hat x_n=1$ and
$\hat x_i$ is the $i$th smallest element of $\{x_1,\ldots,x_{n-1}\}$.
Equip $\N\times (0,1)^\N\times \R^{\N_0}$
with the product of the geometric distribution with parameter $\frac 12$ on
the first coordinate; independent identically
distributed Lebesgue measure on the $x$ coordinates; and independent
 standard normals on the $y$ coordinates, and let
$\mu_\text{PL}$ denote the push-forward measure under the map $\Phi$ on the space $C[0,1]$.

\begin{thm}
Let $\mu_\text{PL}$ be as above. For $\mu_\text{PL}^\N$-a.e.\ subset $V$ of $C[0,1]$,
and each $0<p<1$, almost every realization of $\LARG(V,p)$ is isomorphic to $\mathsf{G}_{TD}$.
\end{thm}

\begin{proof}
Let a randomly sampled $V$ be $\{f_1,f_2,\ldots\}$. We show that
$V$ is transversely dense, so that the desired
result follows from Theorem~\ref{thm:TD}. First notice that the slope of
each piece of a function sampled from $\mu_\text{PL}$
is continuously distributed, so that almost surely, no two pieces of
any $f_i$ and $f_j$ have the same slope. This ensures that
$\cross(f_i,f_j)$ is finite for each $i,j$. As before, the distributions
of $f(0)$ and $f(1)$ are continuous, so with probability
1, no crossings occur at 0 or 1. For any two piecewise linear functions
$f_i$ and $f_j$ with different slopes, local maxima occur
only at the endpoints of the intervals. Given the function $f_j$ and
the $x$ values (but not the $y$ values) defining $f_i$,
the conditional probability that a crossing occurs at one of these values
is 0, so that almost surely no crossings occur at
local maxima or minima of $f_i-f_j$, establishing \ref{it:Tr1}.

To show \ref{it:Tr2}, condition on $f_i$ and $f_j$. These functions almost
surely have finitely many crossings.
Since at each of these crossings, the values of $f_k$ and $f_l$ have
independent continuous distributions, the probability
that $f_k$ crosses $f_i$ and $f_j$ at one of these points; or that
$f_k$ and $f_l$ cross at one of these points is 0.

Finally, let $\mathcal F=\{f_1,\ldots,f_m\}$ and let
$\mathcal F\cup\{f\}$ be transverse. There is a piecewise
linear function $g$ with $\|f-g\|$ arbitrarily small
such that $\CD(f_1,\ldots,f_m,g)=\CD(f_1,\ldots,f_m,f)$. Any suitably close approximation with the
same number of pieces will satisfy $\CD(f_1,\ldots,f_m,h)\sim_\epsilon\CD(f_1,\ldots,f_m,f)$, so that
\ref{it:TD2} follows.
\end{proof}

We exhibit integrated Brownian motion as a third example. Let $\nu$ be the Wiener measure on
standard Brownian motions on $[0,1]$. We then define a map $\Phi\colon \R\times C[0,1]\to C[0,1]$
by
$$
[\Phi(a,f)](x)=a+\int_0^x f(t)\,dt.
$$
We then define a measure $\mu_\text{IBM}$ on $C[0,1]$ that is the push-forward of the product of a
standard normal distribution and $\nu$ under the map $\Phi$. In other words, a function sampled from
$\mu_\text{IBM}$ is obtained by adding a normally-distributed initial condition to the anti-derivative
of a standard Brownian motion.

\begin{thm}
Let $\mu_\text{IBM}$ be as above. For $\mu_\text{IBM}^\N$-a.e.\ subset $V$ of $C[0,1]$,
and each $0<p<1$, almost every realization of $\LARG(V,p)$ is isomorphic to $\mathsf{G}_{TD}$.
\end{thm}

We give a brief proof, primarily emphasizing the parts that differ from the previous examples.

\begin{proof}
Recall that the difference of two independent standard normal random variables is a normal random
variable with variance 2. Similarly, the difference of two independent normal standard Brownian motions
is a Brownian motion with parameters $\mu=0$ and $\sigma^2=2$. Hence, the difference of two
independent realizations of $\mu_\text{IBM}$ is, up to scaling by a factor of $\sqrt 2$, of the same form.
Fubini's theorem implies that the zero set of a Brownian motion almost surely has measure 0, so that
the measure of the set of critical points of a realization of $\mu_\text{IBM}$ is almost surely 0.
Since realizations of Brownian motion on [0,1] are continuous, they are bounded, so that each realization of
$\mu_\text{IBM}$ is almost surely Lipschitz. Conditioned on the Wiener process $f$, it follows that for
all values of $a$ except for a set of measure 0, $\Phi(a,f)$ does not take an integer value at any
critical point. It follows that almost surely, if $f_1$ and $f_2$ are realizations of $\mu_\text{IBM}$,
their difference has non-zero derivative at each element of $\cross(f_1,f_2)$. This ensures
that $\cross(f_1,f_2)$ is closed, and almost surely discrete, hence finite. The remainder of the proof
that a countable set of independently chosen realizations of integrated Brownian motions are
transverse is similar to the previous proofs.

Next, to show that a countable set of independently chosen realizations of integrated Brownian
motions satisfies \ref{it:TD2}, let $\mathcal F=\{f_1,\ldots,f_m\}$ and suppose
$f\in C[0,1]\setminus\mathcal F$ is such that $\mathcal F\cup\{f\}$ is transverse.
Let $\epsilon>0$. There is a $C^1$ function $g$ such that $g'(0)=0$;
$\|g-f\|_{C[0,1]}<\frac\epsilon2$; $\cross(\mathcal F\cup\{g\})=
\cross(\mathcal F\cup\{f\})$; and so that $(f_i-g)'$ is non-zero at each crossing. Then
there is a $C^1$-neighbourhood, $N$, of $g$ so that all elements $h$ of $N$
satisfy $\|f-h\|_{C_0[0,1]}<\epsilon$ and $\CD(f_1,\ldots,f_m,h)\sim_\epsilon
\CD(f_1,\ldots,f_m,f)$. Since the measure on Brownian motions gives positive measure
to every open subset of $\{k\in C_0[0,1]\colon k(0)=0\}$,
the measure of the pairs $(a,f)$ such that $\Phi(a,f)$ lies in $N$
is positive, so that there are almost surely infinitely many $f_k$'s satisfying the conditions
of \ref{it:TD2}.
\end{proof}

\section{IC-dense sets and Brownian motion}\label{secbm}

By definition, transversely dense sets must contain functions that cross finitely often and transversely
with other functions in the set. In this section, we establish results for functions in
$C[0,1]$ that exhibit the opposite type of behaviour: pairs of functions cross infinitely
often. We will call sets of such functions \emph{infinite crossing dense} (IC-dense).
An example of such functions is formed by Brownian motions. We will first show that any two
LARG graphs with IC-dense vertex sets are almost surely isomorphic. We will then define a
natural probability measure on the set of Brownian motions, and show that under this measure
almost all countable sets of Brownian motions are IC-dense.

A countable subset $V$ of $C[0,1]$
will be called \emph{IC-dense}  if it has the following properties:
\begin{enumerate}[label=(IC\arabic*)]
\item The set $V$ is dense in $C[0,1]$;\label{IC1:dense}
\item The sets $\cross(f,g)$ are disjoint for distinct pairs $f,g\in V$;
\label{IC2:notriple}
\item For any $f,g\in V$, $\cross(f,g)$ does not contain 0 or 1; \label{IC3:nocrossbdy}
\item For each distinct pair $f,g\in V$, if  $x\in\cross(f,g)$ for
$x\in (0,1)$, then there exists a sequence of points $x_n$ converging to $x$
monotonically such that $f(x_{2n})-g(x_{2n})<f(x)-g(x)$
and $f(x_{2n+1})-g(x_{2n+1})>f(x)-g(x)$ for each $n$.
\label{IC4:infcross}
\end{enumerate}
We comment on \ref{IC4:infcross}. The requirement that the
sequence $(x_n)$ converges to $x$ monotonically is
to ensure that there is an infinite sequence of crossings of $f$ and
$g$ converging to $x$ (whereas if $f$ and $g$
cross transversely, it is straightforward to build a non-monotonic sequence of $x_n$'s with the alternating property).
The reason that we only require that there is an infinite sequence of $x_n$'s
on one side of $x$ is that the set of crossings of two Brownian motions is perfect
(each crossing is a limit of other crossings), but while some crossings may be
approached by crossings on either side, there are crossings that may be only approached on one side.
To see this, notice that for any pair of
continuous functions that cross, the set of crossings is closed, and therefore has a maximum.

\subsection{IC-dense sets are geometrically Rado}

We will now show that IC-dense sets have the geometric Rado property.
As in the case for transversely dense sets,
the proof is based on the inductive construction of a graph isomorphism. While in
the transversely dense case we required that the isomorphism exactly preserves the
crossing behaviour of the functions, here we wish to preserve the decomposition of the
interval $[0,1]$ into subintervals where pairs of functions cross infinitely often. Before
we state the main theorem, we give definitions and useful lemmas.

{\sl Crossing partitions.} If $\mathcal{F}$ is a finite subset of an IC-dense of $C[0,1]$, then
we build a finite partition of $[0,1]$ into open and closed
sub-intervals by the following procedure. Let
\begin{equation}\label{def:a1}
a_{1}=\min\{x\in[0,1]\colon g(x)-h(x)\in\Z
\text{ for some distinct }
g,h\in\mathcal F\}.
\end{equation}
Suppose $a_1,\ldots,a_n$ have been found, and let $(g_i,h_i,m_i)_{i=1}^n$ be the tuples
in $\mathcal F\times \mathcal F\times \Z$ with $g_i\ne h_i$ so that
$g_i(a_i)-h_i(a_i)=m_i$ for each $1\leq i\leq n$.
We then define
\begin{equation}
\begin{split}
a_{n+1}=\min\big\{ x\in (a_n,1]\colon &
h(x)-g(x)=m\in\Z
\text{ for some distinct }g,h\in{\mathcal F}\\
&\text{such that }\{g,h,m\}\ne\{g_n,h_n,m_n\}\big\}.
\end{split}
\label{def:an}
\end{equation}
If $a_{n+1}$ is defined (that is, the set is not empty), then set
$g_{n+1}=g$, $h_{n+1}=h$ and $m_{n+1}=m$. If not, then the procedure terminates with $a_n$.

Note that it is possible that $\{g_{n+1},h_{n+1}\}=\{g_n,h_n\}$, but
only if they differ by one more or less at $a_{n+1}$ than $a_{n}$. This, in turn,
is only possible if there are exactly two functions. If there are any more functions, then
$g_n$ or $h_n$ would have to cross another element of $\mathcal F$ before crossing each other
again.

\begin{lem}
The procedure to generate the sequence $(a_n)$ as defined in (\ref{def:a1})
and (\ref{def:an}) terminates in a finite number of steps.
\end{lem}

\begin{proof} Suppose, for a
contradiction, that there are infinitely many $a_n$'s. Since $(a_n)$ is
an increasing sequence lying in $[0,1]$, let its limit be $a$. By the
Pigeonhole Principle, there are functions $f^{1}$, $f^{2}$, $f^{3}$ and
$f^{4}$ in $\mathcal F$ (not necessarily distinct, but such that $f^1\ne f^2$,
$f^3\ne f^4$ and $\{f_1,f_2\}\ne \{f_3,f_4\}$) and a
subsequence $(n_j)$ such that $\{g_{n_j},h_{n_j}\}=\{f^{1},f^{2}\}$ and
$\{g_{n_j+1},h_{n_j+1}\}=\{ f^{3},f^{4}\}$ for all $j$. We have $f^{1}(a_{n_j})-f^{2}(a_{n_j})\in\Z$ and
$f^{3}(a_{n_j+1})-f^{4}(a_{n_j+1})\in\Z$ for each $j$.
Taking a limit, and using
continuity of the functions involved, we see that $f^{1}(a)-f^{2}(a)\in\Z$
and $f^{3}(a)-f^{4}(a)\in\Z$. That is, $a\in\cross(f^{1},f^{2})
\cap \cross(f^{3},f^{4})$. This contradicts \ref{IC2:notriple}, establishing
the above procedure terminates as required.
\end{proof}

Now for each $i$ such that $a_i$ is defined, let $b_i=\max
\{x\in [a_i,a_{i+1}]\colon g_i(x)-h_i(x)\in\Z\}$. Note that, since $\cross
(g_i,h_i)$ is closed, this maximum is defined, and since $h_i$ and $g_i$ cross at $b_i$ and a
different pair $h_{i+1},g_{i+1}$ crosses at $a_{i+1}$, we have that $b_i<a_{i+1}$.
The
\emph{crossing partition} $\mathcal P(\mathcal F)$ is then given by
$$
\mathcal P(\mathcal F)=
\{[0,a_1),[a_1,b_1],(b_1,a_2),[a_2,b_2],\ldots,[a_n,b_n],(b_n,1]\},
$$
so that on the \emph{stable intervals}, that is the (relatively) open sets $[0,a_1)$, $(b_1,a_2)$, \ldots,
$(b_{n-1},a_n)$ and $(b_n,1]$, there are no crossings, while on the \emph{crossing intervals}, $[a_i,b_i]$,
$g_i$ and $h_i$ cross infinitely often, but there are no other crossings.

The \emph{crossing data} for $\mathcal F$ consists of the crossing partition together with a
collection of integers $d^I_{f,g}$, one for each pair of distinct elements of $\mathcal F$ and interval
$I\in\mathcal P(\mathcal F)$ except that $d^{[a_i,b_i]}_{g_i,h_i}$
is not defined. Except in this case, for all $x$ in $I$
and all $f,g\in \mathcal F$, then $g(x)-f(x)\in (d^I_{f,g},d^I_{f,g}+1)$. We denote this data by
\cdinf$(\mathcal F)$.
As a consequence, setting $\bar \F=\{f+n\colon f\in\F, n\in\Z\}$,
the elements of $\bar \F$ maintain a constant ordering
on the stable sets. The elements of $\overline{\mathcal F\setminus\{g_i\}}$ and
$\overline{\mathcal F\setminus\{h_i\}}$ maintain a constant
ordering on the crossing intervals.

We write $\cdinf(f_1,\ldots,f_k)\sim_\epsilon\cdinf(g_1,\ldots,g_k)$
if $\mathcal P(\{f_1,\ldots,f_k\})$ contains the same number of intervals as
$\mathcal P(\{g_1,\ldots,g_k\})$; the functions crossing in each crossing interval of
$\mathcal P(\{f_1,\ldots,f_k\})$ have the same indices as the functions crossing in the corresponding
crossing interval of $\mathcal P(\{g_1,\ldots,g_k\})$; $|a'_i-a_i|<\epsilon$
and $|b'_i-b_i|<\epsilon$ for each $i$;
and $d^{I}_{f_i,f_j}=d^{I'}_{g_i,g_j}$ whenever these are defined (where $I$ and $I'$ are corresponding
intervals of $\mathcal P(\{f_1,\ldots,f_k\})$ and $\mathcal P(g_1,\ldots,g_k\})$.
As before $\cdinf(f_1,\ldots,f_k)\sim_1\cdinf(g_1,\ldots,g_k)$ requires that the partitions
have the same number of elements; the same pairs cross in the corresponding crossing intervals
and the gap between a pair of functions has the same floor in each interval of the crossing partition as the
corresponding pair has in the corresponding interval (with the exception of the pair that is crossing in
a crossing interval).

\begin{thm}\label{thm:BMRado}
Let $V$ and $W$ be countable IC-dense subsets of $C[0,1]$ and let $0<p,p'<1$. If $G$ and
$G'$ are sampled from $\LARG(V,p)$ and $\LARG(W,p'),$ respectively, then
$G$ and $G'$ are almost surely isomorphic.

In particular, there exists a graph $\mathsf G_\text{ICD}$ such that if $\mu$ is a measure
on $C[0,1]$ such that $\mu^\N$-a.e. sequence is IC-dense, then $\mu$ is universally geometrically
Rado, and for $\mu^\N$-a.e. vertex set $V$, almost every element of $\LARG(V,p)$ is
isomorphic to $\mathsf G_\text{ICD}$ for any $0<p<1$.
\end{thm}

\begin{proof}
Let $V$ and $W$ be infinite crossing dense subsets of $C([0,1])$ and let $p,p'\in(0,1)$.
We show that a realization of $\LARG(V,p)$ is almost surely isomorphic to a realization of $\LARG(W,p')$.
As in the proof of Theorem~\ref{thm:TD} for transversely dense sets, we construct an isomorphism
via a back-and-forth argument.

The key step of the proof is to show that if $\mathcal F$ and $\mathcal G$ are finite subsets of $V$ and $W$
such that $\cdinf(\mathcal F)\sim_1\cdinf(\mathcal G)$ then given any additional
$f\in V\setminus\mathcal F$,
there is almost surely a $g\in W\setminus\mathcal G$ such that
$\cdinf(\mathcal F\cup \{f\})\sim_1\cdinf(\mathcal G\cup
\{g\})$ and $g$ has the same adjacencies to $\mathcal F$ as $g$ has to $\mathcal G$.

Let $G$ and $G'$ be sampled from $\LARG(V,p)$ and $\LARG(W,p'),$ respectively.
As before, for the purposes of maintaining independence when applying the Borel-Cantelli lemma,
we maintain queues of functions to test, ensuring that we never test whether a pair of vertices is adjacent
more than once.
Let $Q=(f_1,f_2,\ldots)$ be an enumeration of $\mathcal F$ and $Q'=(g_1,g_2,\ldots)$ be an
enumeration of $\mathcal G$.
As before, we build sequences $(i_n)$ and $(j_n)$ satisfying the following:
\begin{enumerate}[label=(T\ref{thm:BMRado}-\arabic*)]
\item
For each $n$, the graphs induced by $G$ on $\{f_{i_1},\ldots,f_{i_n}\}$ and by
$G'$ on $\{g_{j_1},\ldots,g_{j_n}\}$ are isomorphic by mapping $f_{i_l}$ to $g_{j_l}$;
\label{it:infind1}
\item For each $n\in\N$, $\{i_1,\ldots,i_n\}$ and $\{j_1,\ldots,j_n\}$ contain the set
$\{m\in\N\colon m\le \frac n2\}$;
\label{it:infind2}
\item $\cdinf(f_{i_1},\ldots,f_{i_n})\sim_1\cdinf(g_{j_1},\ldots,g_{j_n})$;\label{it:infind3}
\item the remaining elements of the queues $Q$ and $Q'$ at the $n$th stage consist
of all but finitely many elements of $V$ and $W$
respectively. \label{it:infind4}
\end{enumerate}

When $n=1$, we let $i_1=j_1=1$ and remove $f_1$ and $g_1$ from the respective queues.
We focus on a single step in one direction. Suppose that the induction hypotheses are
satisfied after the $(n-1)$st step.
If $n$ is even, let $i_n$ be $\min(\N\setminus\{i_1,\ldots,i_{n-1}\})$, while if $n$ is odd,
let $j_n=\min(\N\setminus\{j_1,\ldots,j_{n-1}\})$. For even $n$, we then seek a suitable $j_n$ to identify
with $i_n$, while if $n$ is odd, we seek a suitable $i_n$. We focus on the case of even $n$,
the case of odd $n$ being essentially identical.

Let $\mathcal F=\{f_{i_1},\ldots,f_{i_{n-1}}\}$ and
$\mathcal G=\{g_{j_1},\ldots,g_{j_{n-1}}\}$. Let the crossing partitions
of $[0,1]$ of $\mathcal F$ and
$\mathcal G$ be
\begin{align*}
\mathcal P(\mathcal F)&=\{[0,a_1),[a_1,b_1],(b_1,a_2),
[a_2,b_2],\ldots,[a_N,b_N],(b_N,1]\}\text{ and}\\
\mathcal P(\mathcal G)&=\{[0,a'_1),[a'_1,b'_1],(b'_1,a'_2),
[a'_2,b'_2],\ldots,[a'_N,b'_N],(b'_N,1]\}.
\end{align*}

Notice that $\mathcal P(\mathcal F\cup\{f_{i_n}\})$ is a refinement
of $\mathcal P(\mathcal F)$: the elements of $\mathcal P
(\mathcal F)$ may be sub-divided in $\mathcal P(\mathcal F\cup\{f_{i_n}\})$
due to intersections of $f_{i_n}$
with elements of $\mathcal F$.
For an element $I$ of $\mathcal P(\mathcal F)$ the sub-partition $\mathcal P_I$ of $I$
is obtained by taking those elements of $\mathcal P(\mathcal F\cup \{f_{i_n}\})$ that lie in $I$. (Note that
$\mathcal P_I$ may be the trivial partition of $I$ if $f_{i_n}$ does not cross any
element of $\mathcal F$ in $I$).

If $I=[a_k,b_k]$ is a crossing interval in $\mathcal P(\mathcal F)$, then two functions $f_{i_l}$ and $f_{i_m}$
cross infinitely often in $I$. By construction, no other functions of $\mathcal F$ cross in $I$.
The sub-partition, $\mathcal P_I$ consists of crossing intervals in which $f_{i_l}$
and $f_{i_m}$ cross infinitely often;
stable intervals (with no crossings); and crossing intervals in which $f_{i_n}$ crosses exactly one
element of $\mathcal F$ (possibly $f_{i_l}$ or $f_{i_m}$), again infinitely often.
The situation is illustrated schematically in Figure \ref{fig:cross}.
If $I$ is a stable interval in $\mathcal P(\mathcal F)$, then $\mathcal P_I$ consists of stable intervals,
and intervals in which $f_{i_n}$ crosses exactly one element of $\mathcal F$ infinitely often.
\begin{figure}
\scalebox{0.6}{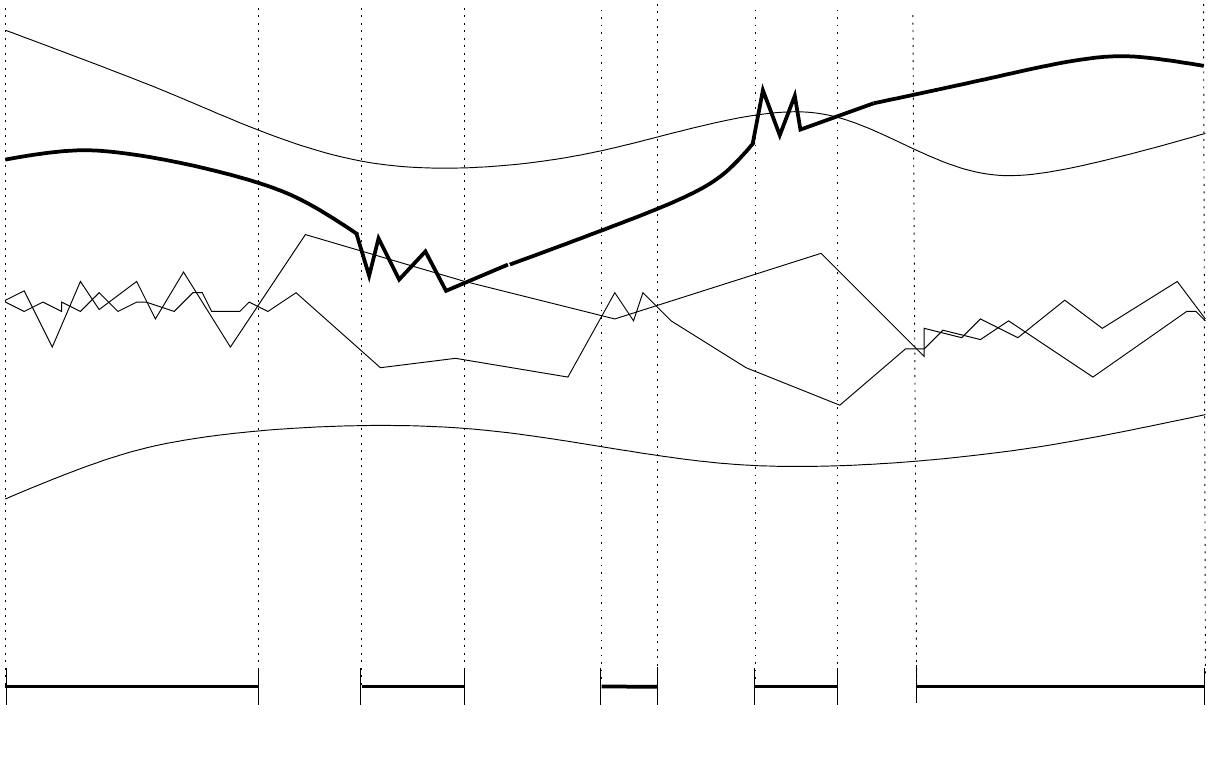}
\caption{Schematic indicating the refinement of a crossing interval $[a_j,b_j]$ in which the
two central functions, $f_{i_l}$ and $f_{i_m}$ are crossing when a new function $f$ (in bold)
is added to the collection. The graph illustrates $f_{i_l,p}$ and $f_{i_m,q}$ for suitable $p,q$
so that the functions agree at the crossing, rather than differing by an integer.
The closed sub-intervals in the refinement are indicated below.}
\label{fig:cross}
\end{figure}

In each crossing interval of $\mathcal P_I$, there are crossings between exactly one pair of functions
in $\mathcal F\cup\{f_{i_n}\}$.
As in the proof of Theorem~\ref{thm:TD}, for $a\in\Z$, we write $f_{i,a}(x)=f_i(x)+a$.

We build a target function $g$ such that all $h\in W\setminus\mathcal G$
sufficiently close to $g$ satisfy $\cdinf(f_{i_1},\ldots,f_{i_n})\sim_1\cdinf(g_{j_1},
\ldots,g_{j_{n-1}},h)$,
which will allow us to maintain \ref{it:infind3}. As in the proof of Theorem~\ref{thm:TD}, this function is built
by interpolating various convex combinations of the elements of $\mathcal G$.

Let $I'$ be the element of $\mathcal P(\mathcal G)$ corresponding to an interval $I$ of
$\mathcal P(\mathcal F)$. We consider separately the cases where
$I$ is a stable interval of $\mathcal P(\mathcal F)$ and where $I$ is a crossing
interval of $\mathcal P(\mathcal F)$.
In the case where $I$ is a stable interval of $\mathcal F$, let $\psi$ be the increasing affine map sending
$I'$ to $I$ and for each $y$, as in Theorem~\ref{thm:TD},
let $\phi_y$ be the piecewise linear function sending $f_{i_k}(\psi(y))+a$ to $g_{j_k}(y)+a$.
The restriction of $g$ to $I'$ is then given by
\begin{equation}\label{eq:defg}
g|_{I'}(y)=\phi_y(f_{i_n}(\psi(y))).
\end{equation}

In the case where $I$ is a crossing interval, suppose that the functions that cross in
$I$ are $f_{i_l}$ and $f_{i_m}$.
We build a target function $g|_{I'}$ in a similar manner to the above, but additional care is needed in
the elements of $\mathcal P_I$ on which $f_{i_l}$ and $f_{i_m}$ cross.

We first build a suitable
partition $\mathcal Q_{I'}$ of $I'$ corresponding to $\mathcal P_I$.
Let the intervals of $\mathcal P_I$ on which $f_{i_l}$ and $f_{i_m}$ cross be
$[\alpha_1,\beta_1],\ldots,[\alpha_s,\beta_s]$ (with $\alpha_1=\min I$ and $\beta_s=\max I$.
On each of the open intervals separating
these crossing intervals, $(\beta_r,\alpha_{r+1})$, $f_{i_m}-f_{i_l}$ takes values in
$(d_r,d_r+1)$ for some $d_r\in\Z$. We then build a partition of $I'$ consisting of sub-intervals
$[\alpha_1',\beta_1'],\ldots, [\alpha'_s,\beta'_s]$ (with $\alpha'_1=\min I'$ and $\beta'_s=\max I'$)
of $I'$ containing all crossings of $g_{j_l}$ and $g_{j_m}$, so that crossings occur at the endpoints
and such that the restriction of $g_{j_m}-g_{j_l}$ to $(\beta'_r,\alpha'_{r+1})$
takes values in $(d_r,d_r+1)$.
This is possible since $\cdinf(f_{i_1},\ldots,f_{i_{n-1}})\sim_1
\cdinf(g_{j_1},\ldots,g_{j_{n-1}})$ and since $W$ satisfies \ref{IC4:infcross}.
Now let $\mathcal Q_{I'}=\{[\alpha'_1,\beta'_1],(\beta'_1,\alpha'_2),\ldots,[\alpha'_s,\beta'_s]\}$
and let $\psi$ be the increasing piecewise affine map that sends each $\alpha'$
and $\beta'$ to the corresponding
$\alpha$ or $\beta$.
Let $J$ be an element of $\mathcal P_I$ and $J'$ be the corresponding element of $\mathcal Q_{I'}$.
If $J$ is a stable interval or an interval
on which $f_{i_n}$ crosses an $f_{i_a}$ other than $f_{i_l}$ or $f_{i_m}$,
then define $\phi_y$ and $g|_{J'}$ as in \eqref{eq:defg}.

If $J$ is an interval on which $f_{i_l}$ and $f_{i_m}$ cross (so that $f_{i_n}$ does not cross either
$f_{i_l}$ or $f_{i_m}$), for $y\in J'$, then let $f_{i_n}-f_{i_l}$ take
values in $(a,a+1)$ and $f_{i_n}-f_{i_m}$ take values in $(b,b+1)$ on $J$.
We then set $u(x)=\max(f_{i_l}(x)+a,f_{i_m}(x)+b)$ and
$v(x)=\min(f_{i_l}(x)+a+1,f_{i_m}(x)+b+1)$ on $J$; and
$u'(y)=\max(g_{j_l}(y)+a,f_{j_m}(y)+b)$ and $v'(y)=\min(g_{j_l}(y)+a+1,f_{j_m}(y)+b+1)$ on $J'$.
Notice that since the only functions crossing in $J$ are $f_{i_l}$ and $f_{i_m}$,
there is a single function, either $u$ or a function of the form $f_{i,c}$ with
$i\in\{1,\ldots,n-1\}\setminus\{l,m\}$ and $c\in\Z$ that is
the closest approximation from below of $f_{i_n}$ throughout the interval $J$.
Similarly, there is a single function
$f_{i,d}$ or $v$ that is the closest approximation from above of $f_{i_n}$.
We call these functions $m(x)$ and $M(x)$.
Let the corresponding functions on $J'$ be $m'(y)$ and $M'(y)$.
Let $\phi_y$ be the affine function sending
$m(\psi(y))$ to $m'(y)$ and $M(\psi(y))$ to $M'(y)$. We then define $g|_{J'}$ as in \eqref{eq:defg}.
One may check that the function $g$ defined in this way is continuous since the maps $\phi_y$ being applied
at the endpoints of the intervals agree. For example in Figure \ref{fig:cross},
$\phi_y$ is defined at $\beta_1$
using $M(x)$ and the function $f$ lying immediately above $f_{i_n}$, while at $\beta_1^+$,
assuming $f_{i_l,a}>f_{i_m,b}$ on $(\beta_1,\alpha_2)$, $\phi_y$ is defined
at $\beta_1^+$ using $f_{i_l,a}$ and $f$.
These agree since $M(\beta_1)=f_{i_l,a}(\beta_1)$.

By construction, in any of the intervals of $\mathcal P(\mathcal G)$,
if $f_{i_n}$ crosses $f_{i_a}$ at $\psi(y)$, then $g$ crosses $g_{j_a}$ at $y$.
Now for any sufficiently small $\epsilon>0$, any $h$ satisfying $\|g-h\|<\epsilon$ has nearby crossings.
For sufficiently small $\epsilon>0$, if $h\in W$ satisfies $\|g-h\|<\epsilon$,
then $\cdinf(f_{i_1},\ldots,f_{i_{n-1}},f_{i_n})\sim_1 \cdinf(g_{j_1},\ldots,g_{j_{n-1}},h)$.
In particular, $\|h-g_{j_a}\|=\|f_{i_n}-f_{i_a}\|$ for each $a$. Each $h$ selected from $Q'$
therefore has probability at least $\min(p',1-p')^{n-1}$ of being adjacent to exactly those $g_{j_a}$'s
corresponding to the $f_{i_a}$'s that are adjacent to $f_{i_n}$. Take the first suitable $g_b$ from the queue
and discard all elements of $Q'$ up to and including $g_b$. Let $j_n=b$ and delete $i_n$ from $Q$.

The induction hypotheses almost surely remain satisfied.
\end{proof}

\subsection{Brownian motion}

We consider the following procedure to select random functions: a
starting point (the value at 0) is chosen from a standard normal
distribution and from there the function follows a standard Brownian
path. Formally, let $(N_n)_{n\in\N}$ be a family of countably many
independent standard normal random variables and let
$(X_n(t))_{n\in\N}$ be a family of countably many independent standard Brownian
motions (with $X_n(0)=0$, $\E X_n(t)=0$ and $\Var X_n(t)=t$ for
each $n\in\N$ and $t\in[0,1]$). Let $f_n(t)=N_n+X_n(t)$ be a random function
obtained in this manner, which we refer to as a {\sl shifted Brownian motion}.

\begin{thm}\label{thm:BM IC-dense}
If $(f_n)_{n\in\N}$ is a family of countably many independent shifted Brownian
motions sampled from the above distribution, then the collection is almost surely IC-dense.
\end{thm}

\begin{proof}
Let $\epsilon>0$ be given. Let $f$ be a continuous function. There exists a
function $g$ where $g(x_i)=f(x_i)$ for $i=0,\ldots,n$ with $0=x_0<x_1<\ldots<x_n=1$
and $g$ is affine in between such that
$\|g-f\|<\frac\epsilon2$. By standard properties of Brownian motion, it is
known that $\PP(|N-y_0|<\frac\epsilon {2n})>0$ and that $\PP\big(|X(x_j)-X(x_{j-1})-(y_j-y_{j-1})|<
\frac\epsilon{2n}\,\big|\,X(x_1),
\ldots,X(x_{j-1})\big)>0$. Finally, it is known that conditioned on the values of
$X(x_1),\ldots,X(x_n)$, the probability that $X$ never strays by more than
$\frac\epsilon2$ from the affine function
joining $X(x_i)$ and $X(x_{i+1})$ is positive (the distribution of
$(X(t))_{x_i\le t\le x_{i+1}}$ conditioned on $X(x_i)$ and $X(x_{i+1})$
is a so-called \emph{Brownian bridge}). Multiplying these
probabilities together, we note that the probability that $\|f_j-f\|<\epsilon$ is
positive, and is the same for each $j$. Hence, by the second Borel-Cantelli lemma,
there almost surely exists a $j$
such that $\|f_j-f\|<\epsilon$. Since $C([0,1])$ is separable, the above argument
shows that $(f_n)$ is almost surely dense.

Condition \ref{IC2:notriple} follows from a corollary of L\'evy's theorem
(Corollary~2.23 in \cite{MortersPeres}) which states that
for any non-zero point $p$ in $\R^2$, the standard 2-dimensional Brownian motion
almost surely does not pass through $p$. To see \ref{IC2:notriple} from this,
we deal with two cases: we first show that there
are no three functions all crossing at a single $x$ value; and then show that
there are not two distinct pairs of functions with both pairs crossing at a single $x$
value. For the first of these, it suffices (since there are countably many triples)
to show that almost surely, $f_1$, $f_2$ and $f_3$ do not
all cross at any $x$. That is, it suffices to show that for any
integers $m$ and $n$, with probability 1, there does not exist an $x$
such that $f_2(x)-f_1(x)=m$ and $f_3(x)-f_1(x)=n$. Conditional on
$N_1=y_1$, $N_2=y_2$ and $N_3=y_3$, the values at $x=0$,
we therefore have to show that
$$
\PP(\exists x\colon X_2(x)-X_1(x)=m+y_2-y_1\text{ and }
X_3(x)-X_1(x)=n+y_3-y_1)=0.
$$
This does not immediately follow from the theorem mentioned above
about 2-dimensional Brownian motion since $((X_2-X_1)(x),(X_3-X_1)(x))$
is not a standard Brownian motion: the coordinates are correlated. It is,
however, expressible as a fixed
linear transformation of a standard Brownian motion (see, for example,
\cite{MortersPeres} page 14) and hence, the previous theorem applies.

Similarly, to show that $f_1$, $f_2$, $f_3$ and $f_4$ do not have the
property that $f_1$ and $f_2$ cross at some $x$ and $f_3$ and $f_4$
cross at the same $x$, we need to show that
$$
\PP(\exists x\colon X_2(x)-X_1(x)=m+y_1-y_2\text{ and }
X_4(x)-X_3(x)=n+y_3-y_4)=0.
$$
This follows directly from the quoted theorem since $((X_2-X_1)(x),
(X_4-X_3)(x))$ is a scaled copy of the standard Brownian motion (by a scale
factor of $\sqrt 2$).

To see that $f_1$ and $f_2$ do not cross at 0, we observe that
$\PP(N_2=N_1+n|N_1=x)=0$ and so integrating over $x$, the
probability of a crossing at 0 is 0. A similar argument involving
$N_1$, $N_2$, $X_1(1)$ and $X_2(1)$ shows that the probability
of a crossing at 1 is 0, and this proves \ref{IC3:nocrossbdy}.

To verify the last condition, notice that $f_2(x)-f_1(x)=X_2(x)-X_1(x)+N_2-N_1$.
By Theorem~2.11 of \cite{MortersPeres}, a Brownian
motion has countably many local extrema and by Theorem~5.14 of \cite{MortersPeres},
no local points of increase or decrease (a point of local increase for a Brownian motion
$Z(x)$ is an $x_0$ such that for some $\delta>0$,
on $(x_0,x_0+\delta)$, $Z(x)>Z(x_0)$ and
on $(x_0-\delta,x_0)$, $Z(x)<Z(x_0)$, with a point of local decrease being
defined similarly).
Given $X_1$ and $X_2$, the
probability that $N_2-N_1$ is such that value taken at one of the countably many local
extrema of $f_2-f_1$ is an integer is 0 (since this only occurs
for countably many values of $N_2-N_1$ and $N_1$ and $N_2$ are independent and have
continuous distributions). Hence, we have that almost surely, each crossing
point of $f_2-f_1$ is not a local extremum or a point of local increase or decrease.
Since these possibilities are ruled out, the remaining possibility is
the IC-density condition \ref{IC4:infcross}.
\end{proof}

\begin{cor}\label{cor:RadoBM}
If $V$ is a countable collection of independent samples of shifted Brownian motion,
then for all $p\in (0,1)$, $\LARG(V,p)$ is almost surely isomorphic to
$\mathsf G_\text{ICD}$.
\end{cor}

\begin{proof}
By Theorem~\ref{thm:BM IC-dense}, a countable collection of realizations of Brownian
motions with random starting points is almost surely IC-dense. By Theorem~\ref{thm:BMRado},
if one builds a LARG graph from such a realization, it is almost surely isomorphic to
$\mathsf G_\text{ICD}$.
\end{proof}

\section{Non-isomorphism}

In the previous sections we showed that LARG graphs on transversely dense sets are
almost surely isomorphic to $\mathsf G_\text{TD}$, while LARG
graphs on IC-dense sets are almost surely isomorphic to  $\mathsf G_\text{ICD}$.
In this section, we will show that these two graphs are not isomorphic.
By Corollary~\ref{cor:noniso}, we already know they are not isomorphic to the infinite
graphs obtained from LARG graphs on other Banach spaces.

Before we give this theorem, we recall some classical results which will be
used in the proof. A theorem of Mazur and Ulam \cite{mu} states that a
surjective isometry of real normed spaces is necessarily affine.
Dilworth \cite{Dilworth} extended this by showing that an approximate
surjective isometry of Banach spaces may be uniformly approximated by
a linear isometry. Recall from Lemma~\ref{lem:wellcon}
that any isomorphism between well-connected LARG graphs
satisfies $\lceil \|x-y\|\rceil=\lceil \theta(x)-\theta(y)\|\rceil$. The
following corollary of Dilworth's result states that any step-isometry
gives rise to a surjective linear isometry that bounds the distance
between pre-image and image.

\begin{thm}\label{thm:Dilworth}(\cite{BJQ1})
Suppose $(f_n)$ and $(g_n)$ are dense sequences in separable Banach
spaces $X$ and $Y$ and $\lceil \|f_i-f_j\|\rceil =
\lceil \|g_i-g_j\|\rceil$ for each $i,j$.
We then have that there exists a surjective
linear isometry, $L$, from $X$ to $Y$ and a constant $\kappa>0$ such that
$\|g_i-L(f_i)\|<\kappa$
for all $i$.
\end{thm}

The following classical result classifies the surjective isometries between Banach
spaces of real continuous functions.

\begin{thm}[Banach--Stone]\label{thm:BS}
If $A$ and $B$ are compact Hausdorff spaces and $L$ is a surjective linear isometry from
$C(A,\R)$
to $C(B,\R)$, then there exists a homeomorphism $\varphi\colon B\to A$ and a function
$g\in C(B,\R)$ with $|g(y)|=1$ for all $y$ such that
$$
L[f](y)=g(y)f(\varphi(y))\text{ for all $f\in C(A,\R)$ and $y\in B$.}
$$
\end{thm}

We apply this in the special case that $X$ and $Y$ are both the closed unit interval,
so that $g$ is a constant function, either $g\equiv1$ or $g\equiv-1$.
This will lead to necessary conditions for two sets to give rise to isomorphic
graphs, as stated in the following theorem.

\begin{thm}\label{thm:radovrado}
If $\{f_n\colon n\in\N\}$ and $\{g_n\colon n\in\N\}$ are countable dense sets in
$C([0,1])$; and $G_1$ and $G_2$ are well-connected geometric graphs
with these vertex sets, then there exists a homeomorphism $\phi$ from $[0,1]$
to itself and a permutation of $\N$, $i\mapsto n_i$, such that one of the following holds:
\begin{enumerate}
\item
(Order-preserving case)
For each $x_0\in [0,1]$,
if $f_{i}(x_0)> f_{j}(x_0)$ then $g_{n_i}(\phi(x_0))\ge g_{n_j}(\phi(x_0))$;
\item
(Order-reversing case)
For each $x_0\in [0,1]$, if $f_{i}(x_0)> f_j(x_0)$ then $g_{n_i}(\phi(x_0))\le g_{n_j}(\phi(x_0))$.
\end{enumerate}
\end{thm}

The idea of the proof (in the order-preserving case)
is that we find a function $f_k$ where $f_i-f_k$ and $f_j-f_k$ are non-negative,
$\lfloor \|f_i-f_k\|\rfloor>\lfloor \|f_j-f_n\|\rfloor$ and $f_i-f_k$ and $f_j-f_k$ take values close to their
maximum only on a small neighbourhood of $x_0$. We use the facts that
$g_a\approx f_a\circ\psi$ and $\lfloor \|g_{n_a}-g_{n_b}\|\rfloor=\lfloor \|f_a-f_b\|
\rfloor$ for each $a,b$ to deduce that $g_{n_i}$ dominates $g_{n_j}$ somewhere
arbitrarily close to $\psi^{-1}(x_0)$.

\begin{proof}
Let $\{f_n\}$ and $\{g_n\}$ be as described and let $G_1$ and $G_2$ be isomorphic
well-connected geometric graphs on these sets. We may re-label $(g_n)$ so that for each $n$,
$f_n$ is identified with $g_n$ in the isomorphism. By Lemma~\ref{lem:wellcon},
such an isomorphism must almost surely satisfy:
$\lceil \|f_i-f_j\|\rceil = \lceil \|g_i-g_j\|\rceil$
for each $i,j$. By Theorems \ref{thm:Dilworth} and \ref{thm:BS},
there exist a self-homeomorphism $\psi$ of $[0,1]$, a constant
$\rho=\pm 1$ and a $\kappa>0$ such that
\begin{equation}
\|\rho f_i\circ\psi-g_i\|<\kappa\text{ for all $i$.}\label{eq:BScor}
\end{equation}
We deal with the order-preserving case, $\rho=1$, although the case $\rho=-1$
is analogous. Let $\phi=\psi^{-1}$.

Let $x_0\in [0,1]$ and suppose that $f_i(x_0)>f_j(x_0)$. We let
$\Delta=f_i(x_0)-f_j(x_0)$. For any $l,N\in\N$, define
$$
h_{l,N}(x)=\max\big(f_i(x),f_j(x)+\Delta\big)-N-\tfrac\Delta 2+l|x-x_0|.
$$
Given $l\in\N$, let $N_l\in\N_{\ge 2}$ be chosen
sufficiently large that
\begin{equation}\label{eq:hbound}
h_{l,N_l}(x)<\min\big(f_i(x),f_j(x)\big)-2\kappa-\Delta\text{ for all $x\in[0,1]$.}
\end{equation}
We write $h$ for $h_{l,N_l}$. We now have that, for each $x\in [0,1]$,
\begin{align*}
2\kappa+\Delta&\le f_i(x)-h(x)\le N+\tfrac\Delta 2-l|x-x_0|\\
2\kappa+\Delta&\le f_j(x)-h(x)\le N-\tfrac\Delta 2-l|x-x_0|,
\end{align*}
where the upper bounds for $f_i-h$ and $f_j-h$ are equalities at $x=x_0$.
In particular, we see $\|f_i-h\|=N+\frac\Delta2$ and $\|f_j-h\|=N-\frac\Delta2$, so
$\lfloor\| f_i-h\|\rfloor\ge  N$ while $\|\lfloor f_j-h\|\rfloor < N$.

Let $k\in\N$ be such that $\|f_k-h\|<\frac\Delta2$. This implies
that $\lfloor \|f_i-f_k\|\rfloor\ge  N$
and $\lfloor \|f_j-f_k\|\rfloor< N$. By Lemma~\ref{lem:wellcon}, we
have $\lfloor \|g_i-g_k\|\rfloor\ge N$ and $\lfloor \|g_j-g_k\|\rfloor< N$ also.

By \eqref{eq:BScor}, we have that
\begin{align*}
	&-\kappa<f_n\circ\psi(y)-g_n(y)<\kappa\text{ for $n\in\N$ and $y\in[0,1]$; so}\\
	&-\kappa-\tfrac\Delta2<h\circ\psi(y)-g_k(y)<\kappa+\tfrac\Delta2\text{ for $y\in[0,1]$}.
\end{align*}
Combining this with \eqref{eq:hbound}, it follows that $g_i-g_k$ and $g_j-g_k$ are non-negative.

We also see that
\begin{align*}
g_i(y)-g_k(y)&\le f_i(\psi(y))-h(\psi(y))+2\kappa+\tfrac\Delta 2\\
&\le N+\Delta+2\kappa-l|\psi(y)-x_0|
\end{align*}
Notice that $0\le g_i(y)-g_k(y)< N$ whenever $|\psi(y)-x_0|> (\Delta+2\kappa)/l$.
Since $\|g_i-g_k\|\ge N$, it follows that there exists $y_l$ such that
$|\psi(y_l)-x_0|\le (\Delta+2\kappa)/l$ with
$g_i(y_l)-g_k(y_l)\ge N$. On the other hand, since $\|g_j-g_k\|<N$, we see
$g_j(y_l)-g_k(y_l)< N$. In particular, we deduce that $g_i(y_l)>g_j(y_l)$.
We now have a sequence of points $y_l$ where
$g_i(y_l)>g_j(y_l)$. We also have $|\psi(y_l)-x_0|\to 0$ as $l\to\infty$,
which implies $y_l\to y_0=\phi(x_0)$. By continuity of $g_i$ and $g_j$,
we deduce that $g_i(\phi(x_0))\ge g_j(\phi(x_0))$.
\end{proof}

\begin{cor}\label{BMnonpoly}
The  graphs $\mathsf{G}_\text{ICD}$ and $\mathsf G_\text{TD}$ are non-isomorphic.
\end{cor}

\begin{proof}
Let $(f_n)$ be IC-dense and $(g_n)$ be a transversely dense sequence of polynomials
(such a sequence exists by Theorem~\ref{thm:polytd}).
Let $G_1$ and $G_2$ be LARG graphs with vertices $(f_n)$ and
$(g_n)$, respectively, so that $G_1$ is almost surely isomorphic to
$\mathsf{G}_\text{ICD}$ and $G_2$ is almost surely isomorphic to $\mathsf{G}_\text{TD}$.
Suppose for a contradiction that $G_1$ and $G_2$ are isomorphic.

We now apply Theorem~\ref{thm:radovrado}, and assume without loss of generality
that we are in the orientation-preserving case. Thus,
there exists a permutation $i\mapsto n_i$ of $\N$ and a homeomorphism $\psi$ of
$[0,1]$ such that $f_i(x)>f_j(x)$ implies $g_{n_i}(\psi(x)) \ge g_{n_j}(\psi(x))$ for
all $x$. By the infinite crossing dense condition, there must be two functions
that intersect. That is, there exist $i\ne j$ so  that $f_i(x)=f_j(x)$ for some
$x\in (0,1)$. By \ref{IC4:infcross}, there exists a strictly monotonic sequence $x_m\to x$ such that
$f_i(x_m)=f_j(x_m)$ and between any pair $x_m$ and $x_{m+1}$, there exist
$y_m$ and $z_{m}$ such that $f_i(y_m)>f_j(y_m)$ and $f_i(z_m)<f_j(z_m)$.
It follows that $g_{n_i}(\psi(y_m))\ge
g_{n_j}(\psi(y_m))$ and $g_{n_i}(\psi(z_m))\ge g_{n_j}(\psi(z_m))$.
In particular, by the Intermediate Value Theorem, there exists $w_m$
between $y_m$ and $z_m$ such that
$g_{n_i}(\psi(w_m))=g_{n_j}(\psi(w_m))$. Hence, the two polynomials
$g_{n_i}$ and $g_{n_j}$ agree at infinitely many points, and so are equal,
which is a contradiction.
\end{proof}



%


\section*{Acknowledgments} 
We thank the referee for a careful reading of the paper, and for making helpful suggestions
regarding the presentation.

%
%

\begin{dajauthors}
\begin{authorinfo}[ab]
  Anthony Bonato\\
  Department of Mathematics\\
  Ryerson University\\
  Toronto, ON, CANADA\\
  abonato\imageat{}ryerson\imagedot{}ca \\
\end{authorinfo}
\begin{authorinfo}[jj]
  Jeannette Janssen\\
  Department of Mathematics and Statistics\\
  Dalhousie University\\
  Halifax, NS, CANADA\\
  jeannette\imagedot{}janssen\imageat{}dal\imagedot{}ca \\
\end{authorinfo}
\begin{authorinfo}[aq]
  Anthony Quas\\
  Department of Mathematics and Statistics\\
  University of Victoria\\
  Victoria, BC, CANADA\\
  aquas\imageat{}uvic\imagedot{}ca\\
\end{authorinfo}
\end{dajauthors}

\end{document}